\documentclass[10pt]{article}
\usepackage{latexsym}
\usepackage{amsfonts}
\usepackage{enumerate}
\usepackage{multicol}
\usepackage{graphicx}
\usepackage{amssymb}
\usepackage{amsmath}
\usepackage{epic}

\title{The Maximum Size of $(k,l)$-Sum-Free Sets in Cyclic Groups \\[.4in]}

\author{B\'{e}la Bajnok\footnote{Corresponding author} \\[.1in] {\small Department of Mathematics, Gettysburg College, U.S.A.} \\
{\small E-mail:  bbajnok@gettysburg.edu} \\ [.2in]
and \\[.2in]
Ryan Matzke \\[.1in] {\small School of Mathematics, University of Minnesota, U.S.A.} \\
{\small E-mail:  matzk053@umn.edu} \\ [.2in]
 \\[.4in]}

\date{August 21, 2018}

\newtheorem{thm}{Theorem}

\newtheorem{lem}[thm]{Lemma}
\newtheorem{cor}[thm]{Corollary}

\begin{document}

\maketitle

\begin{abstract}

A subset $A$ of a finite abelian group $G$ is called $(k,l)$-sum-free if the sum of $k$ (not-necessarily-distinct) elements of $A$ never equals the sum of $l$ (not-necessarily-distinct) elements of $A$.  We find an explicit formula for the maximum size of a $(k,l)$-sum-free subset in $G$ for all $k$ and $l$ in the case when $G$ is cyclic by proving that it suffices to consider $(k,l)$-sum-free intervals in subgroups of $G$.  This simplifies and extends earlier results by Hamidoune and Plagne and by Bajnok.

\end{abstract}



\thispagestyle{empty}

\section{Introduction}

Let $G$ be an additively written abelian group of finite order $n$ and exponent $e(G)$.   When $G$ is cyclic, we identify it with $\mathbb{Z}_n=\mathbb{Z}/n\mathbb{Z}$; we consider $0,1,\dots,n-1$ interchangeably as integers and as elements of $\mathbb{Z}_n$.  

For subsets $A$ and $B$ of $G$, we use the standard notations of $A+B$ and $A-B$ to denote the set of two-term sums and differences, respectively, with one term chosen from $A$ and the other from $B$.  If, say, $A$ consists of a single element $a$, we simply write $a+B$ and $a-B$ instead of $A+B$ and $A-B$.  For a subset $A$ of $G$ and a positive integer $h$, $hA$ denotes the $h$-fold {\em sumset} of $A$, that is, the collection of $h$-term sums with (not-necessarily-distinct) elements from $A$.  Note that the $h$-fold sumset of $A$ is (usually) different from its $h$-fold {\em dilation} $h \cdot A =\{ha \mid a \in A\}$.

For positive integers $k$ and $l$, with $k>l$, we call a subset $A$ of $G$ {\em $(k,l)$-sum-free} if $kA$ and $lA$ are disjoint or, equivalently, if $$0 \not \in kA-lA.$$  For example, we see that $A=\{1,2\}$ is a $(5,2)$-sum-free set in $\mathbb{Z}_9$: We have $5A=\{5,6,7,8,0,1\}$ and $2A=\{2,3,4\}$.  (In this example, $kA$ and $lA$ are not only disjoint, but also partition the group; such $(k,l)$-sum-free sets are called {\em complete}.)  We denote the maximum size of $(k,l)$-sum-free subsets in $G$ by $\mu (G, \{k,l\})$.     As our main result in this paper, we determine $\mu (\mathbb{Z}_n, \{k,l\})$ for all $n$, $k$, and $l$.  

Before we state our results, it may be interesting to briefly review the history of this problem.  A $(2,1)$-sum-free set is simply called a {\em sum-free} set.  Sum-free sets in abelian groups were first introduced by Erd\H{o}s in \cite{Erd:1965a} and then studied systematically by Wallis, Street, and Wallis in \cite{WalStrWal:1972a}.   

We can construct sum-free sets in $G$ by selecting a subgroup $H$ in $G$ for which $G/H$ is cyclic and then taking the ``middle one-third'' of the cosets of $H$.  More precisely, with $d$ denoting the index of $H$ in $G$, we see that 
\begin{eqnarray*}
A & = & \bigcup_{i=\lceil (d-1)/3 \rceil}^{2\lceil (d-1)/3 \rceil -1} (i+H)
\end{eqnarray*} is sum-free in $G$, and thus
\begin{eqnarray*} 
\mu (G, \{2,1\}) &\geq & \max_{d |e(G)} \left \{ \left \lceil \frac{d-1}{3} \right \rceil \cdot \frac{n}{d}  \right \}.
\end{eqnarray*}    Using a version of Kneser's Theorem, Diamanda and Yap proved that we cannot do better in cyclic groups:

\begin{thm} [Diamanda and Yap, 1969; cf.~\cite{DiaYap:1969a}, \cite{WalStrWal:1972a}]
For all positive integers $n$, we have 
\begin{eqnarray*} 
\mu (\mathbb{Z}_n, \{2,1\}) &= & \max_{d |n} \left \{ \left \lceil \frac{d-1}{3} \right \rceil \cdot \frac{n}{d}  \right \}.
\end{eqnarray*}
\end{thm}

The fact that the lower bound is also exact in the case of noncyclic groups was established by Green and Ruzsa via complicated methods that, in part, also relied on a computer:

\begin{thm} [Green and Ruzsa, 2005; cf.~\cite{GreRuz:2005a}]  \label{Green and Ruzsa thm}
For any abelian group $G$ of order $n$ and exponent $e(G)$, we have 
\begin{eqnarray*} 
\mu (G, \{2,1\}) &= & \max_{d |e(G)} \left \{ \left \lceil \frac{d-1}{3} \right \rceil \cdot \frac{n}{d}  \right \}.
\end{eqnarray*}
\end{thm}

The first result for general $k$ and $l$ was given by Bier and Chin:

\begin{thm} [Bier and Chin, 2001; cf.~\cite{BieChi:2001a}]
Let $p$ be a positive prime.  If $k-l$ is divisible by $p$, then $\mu (\mathbb{Z}_p, \{k,l\})=0$, otherwise
\begin{eqnarray*}
\mu (\mathbb{Z}_p, \{k,l\}) &=& \left \lceil \frac{p-1}{k+l} \right \rceil.
\end{eqnarray*} 
\end{thm}

This was generalized by Hamidoune and Plagne:

\begin{thm} [Hamidoune and Plagne, 2003; cf.~\cite{HamPla:2003a}]  \label{mu when n and k-l are rel prime}
If $k-l$ is relatively prime to $n$, then 
 \begin{eqnarray*} 
\mu (\mathbb{Z}_n, \{k,l\}) &=& \max_{d |n} \left \{ \left \lceil \frac{d-1}{k+l} \right \rceil \cdot \frac{n}{d}  \right \}.
\end{eqnarray*}
\end{thm}

The case when $n$ and $k-l$ are not relatively prime is considerably more complicated.  We have the following bounds of the first author:
\begin{thm} [Bajnok, 2009; cf.~\cite{Baj:2009a}]  \label{mu bounds lower upper}
For all positive integers $n$, $k$, and $l$, with $k>l$, we have
 \begin{eqnarray*} 
\max_{d |n} \left \{ \left \lceil \frac{d-\delta}{k+l} \right \rceil \cdot \frac{n}{d}  \right \} & \leq & \mu (\mathbb{Z}_n, \{k,l\}) \; \; \leq \; \; \max_{d |n} \left \{ \left \lceil \frac{d-1}{k+l} \right \rceil \cdot \frac{n}{d}  \right \},
\end{eqnarray*}
where $\delta=\gcd(d,k-l)$.
\end{thm}

Until now, not even a conjecture was known for the actual value of $\mu (\mathbb{Z}_n, \{k,l\})$; as our example above demonstrates, both inequalities in Theorem \ref{mu bounds lower upper} may be strict: $\mu (\mathbb{Z}_9, \{5,2\})=2$, while the lower and upper bounds above are $1$ and $3$, respectively. 
Here we prove the following result:

\begin{thm}  \label{mu all n k l}
For all positive integers $n$, $k$, and $l$, with $k>l$, we have
\begin{eqnarray*} \mu (\mathbb{Z}_n, \{k,l\}) & = & \max_{d |n} \left \{ \left \lceil \frac{d-(\delta-r)}{k+l} \right \rceil \cdot \frac{n}{d}  \right \},
\end{eqnarray*}
where $\delta=\gcd(d,k-l)$ and $r$ is the remainder of $l  \lceil (d-\delta)/(k+l) \rceil$ mod $\delta$.

\end{thm}
We may observe that $\delta -r$ is between 1 and $\delta$, inclusive, so Theorem \ref{mu bounds lower upper} follows from Theorem \ref{mu all n k l}; in particular, we get Theorem \ref{mu when n and k-l are rel prime} when $n$ and $k-l$ are relatively prime.  We also see that Theorem \ref{mu all n k l} implies that $\mu (\mathbb{Z}_n, \{k,l\}) =0$ if, and only if, $k-l$ is divisible by $n$.

Let us now turn to the discussion of our approach.  The main role in our development will be played by {\em arithmetic progressions}, that is, sets of the form
$$A=\{a+ i \cdot b \mid i=0,1,\dots,m-1\}$$ for some positive integer $m$ and elements $a$ and $b$ of $\mathbb{Z}_n$.  (We assume that $m \leq n/\gcd(n,b)$ and thus $A$ has size $|A|=m$.  Note also that $a$ and $b$ are not uniquely determined by $A$; the only time when this will make a difference for us is when $|A|=1$, in which case we set $b=1$.)  In \cite{HamPla:2003a}, Hamidoune and Plagne proved that, if $n$ and $k-l$ are relatively prime, then $\mu (\mathbb{Z}_n, \{k,l\})$ equals     
$$\max_{d |n} \left\{ \alpha (\mathbb{Z}_d, \{k,l\}) \cdot \frac{n}{d}  \right \},$$ where $\alpha (\mathbb{Z}_d, \{k,l\})$ is the maximum size of a $(k,l)$-sum-free arithmetic progression in $\mathbb{Z}_d$.  Hamidoune and Plagne only treat the case when $n$ and $k-l$ are relatively prime; as they write, ``in the absence of this assumption, degenerate behaviors may appear.'' Nevertheless, as the first author proved, the identity remains valid in the general case:

\begin{thm} [Bajnok, 2009; cf.~\cite{Baj:2009a}]  \label{Bajnok alpha}

For all positive integers $n$, $k$, and $l,$ with $k>l$, we have

\begin{eqnarray*} 
\mu (\mathbb{Z}_n, \{k,l\}) & = & \max_{d |n} \left\{ \alpha (\mathbb{Z}_d, \{k,l\}) \cdot \frac{n}{d}  \right \}.
\end{eqnarray*}

\end{thm}

When attempting to evaluate $\alpha (\mathbb{Z}_d, \{k,l\})$, one naturally considers two types of arithmetic progressions: those with a common difference $b$ that is not relatively prime to $d$ (in which case the set is contained in a coset of a proper subgroup), and those where $b$ is relatively prime to $d$ (in which case the set, unless of size $1$, is not contained in a coset of a proper subgroup).  Accordingly, in \cite{HamPla:2003a} Hamidoune and Plagne define $\beta (\mathbb{Z}_d, \{k,l\})$ as the maximum size of a $(k,l)$-sum-free arithmetic progression with $\gcd(b,d)>1$, and $\gamma (\mathbb{Z}_d, \{k,l\})$ as the maximum size of a $(k,l)$-sum-free arithmetic progression with $\gcd(b,d)=1.$  Clearly,
$$\alpha (\mathbb{Z}_d, \{k,l\}) = \max \left \{ \beta (\mathbb{Z}_d, \{k,l\}) , \gamma (\mathbb{Z}_d, \{k,l\}) \right \}.$$

The authors of \cite{HamPla:2003a} evaluate both $\beta (\mathbb{Z}_d, \{k,l\})$ and $\gamma (\mathbb{Z}_d, \{k,l\})$ under the assumption that $d$ and $k-l$ are relatively prime.  We are able to find $\gamma (\mathbb{Z}_d, \{k,l\})$ without this assumption:

\begin{thm} \label{thm gamma eval}

For all positive integers $d$, $k$, and $l$, with $k>l$, we have
\begin{eqnarray*}  \gamma (\mathbb{Z}_d, \{k,l\}) & = &  \left \lceil \frac{d-(\delta-r)}{k+l} \right \rceil , \end{eqnarray*} 
where $\delta=\gcd(d,k-l)$ and $r$ is the remainder of $l  \lceil (d-\delta)/(k+l) \rceil$ mod $\delta$.

\end{thm}

However, evaluating $\beta (\mathbb{Z}_d, \{k,l\})$ in general does not seem feasible.  Luckily, as we here prove, this is not necessary, since we have the following result:

\begin{thm} \label{gamma is enough}

For all positive integers $n$, $k$, and $l$ with $k>l$ we have
\begin{eqnarray*} \max_{d |n} \left\{ \alpha (\mathbb{Z}_d, \{k,l\}) \cdot \frac{n}{d}  \right \} & = &  \max_{d |n} \left\{ \gamma (\mathbb{Z}_d, \{k,l\}) \cdot \frac{n}{d}  \right \}. 
\end{eqnarray*}

\end{thm}

Therefore, Theorem \ref{mu all n k l} follows readily from Theorems \ref{Bajnok alpha}, \ref{gamma is enough}, and \ref{thm gamma eval}.  In Sections \ref{section thm gamma eval} and \ref{section gamma is enough} below we prove Theorems \ref{thm gamma eval} and \ref{gamma is enough}, respectively.  In Section \ref{Further questions} we discuss some further related questions.

\section{The Maximum Size of $(k,l)$-Sum-Free Intervals}  \label{section thm gamma eval}

In this section we evaluate $ \gamma (\mathbb{Z}_d, \{k,l\})$ and thus prove Theorem \ref{thm gamma eval}.  Note that if 
$$A=\{a+ i \cdot b \mid i=0,1,\dots,m-1\},$$
with $b$ relatively prime to $d$, then $b \cdot c=1$ for some $c \in \mathbb{Z}_d$, and thus the $c$-fold dilation
$$c \cdot A = \{ca+ i  \mid i=0,1,\dots,m-1\}$$ of $A$ is the interval $[ca,ca+m-1]$; furthermore, $A$ is $(k,l)$-sum-free in $\mathbb{Z}_d$ if, and only if, $c \cdot A$ is.  Therefore, we may restrict our attention to intervals.

First, we prove a lemma.

\begin{lem}  \label{lem iff k,l free interval}

Suppose that $k$, $l$, and $d$ are positive integers and that $k>l$; let $\delta=\gcd(d,k-l)$.  Then  $\mathbb{Z}_d$ contains a $(k,l)$-sum-free interval of size $m$ if, and only if,
\begin{eqnarray*}  k(m-1) + \left \lceil (l(m-1)+1)/ \delta \right \rceil \cdot \delta & < & d.  \end{eqnarray*}

\end{lem}

In particular, $\mu (\mathbb{Z}_d, \{k,l\})=0$ if, and only if, $d$ divides $k-l$.

{\em Proof:}  Let $A=[a,a+m-1]$ with $a \in \mathbb{Z}_d$ and $|A|=m$.  (As customary, our notation stands for the interval $\{a,a+1,\dots,a+m-1\}$.)  Note that $A$ is $(k,l)$-sum-free if, and only if, $$0 \not \in kA-lA.$$  Observe that $kA-lA$ is also an interval, namely
$$kA-lA=[(k-l)a-l(m-1), \; (k-l)a + k(m-1)].$$
Therefore, $A$ is $(k,l)$-sum-free if, and only if,  there is a positive integer $b$ for which
$$(k-l)a-l(m-1) \geq bd+1$$ and
$$(k-l)a + k(m-1) \leq (b+1)d-1.$$  The set of these two inequalities is equivalent to
$$l(m-1)+1 \leq (k-l)a-bd \leq d-k(m-1)-1,$$ or
$$\frac{l(m-1)+1}{\delta} \leq \frac{(k-l)}{\delta} \cdot a-\frac{d}{\delta} \cdot b  \leq \frac{d-k(m-1)-1}{\delta}.$$
Here $\tfrac{(k-l)}{\delta}$ and $\tfrac{d}{\delta}$ are relatively prime, so every integer can be written in the form $$\frac{(k-l)}{\delta} \cdot a-\frac{d}{\delta} \cdot b$$ for some $a$ and $b$; we may also assume that $0 \leq a \leq d/\delta-1$ and hence $0 \leq a \leq d-1$.  Therefore, $\mathbb{Z}_d$ contains a $(k,l)$-sum-free interval of size $m$ if, and only if, there is an integer $C$ with 
$$\frac{l(m-1)+1}{\delta} \leq C  \leq \frac{d-k(m-1)-1}{\delta},$$ or, equivalently,
$$\left \lceil \frac{l(m-1)+1}{\delta} \right \rceil \leq  \frac{d-k(m-1)-1}{\delta},$$ which is further equivalent to
$$ k(m-1)+ \left \lceil ( l(m-1)+1)/\delta \right \rceil \cdot \delta  < d,$$ as claimed.
$\Box$

{\em Proof of Theorem \ref{thm gamma eval}:}  Let $\gamma_d=\gamma (\mathbb{Z}_d, \{k,l\})$, $$f=\left \lceil \frac{d- \delta}{k+l} \right \rceil,$$ and $$m_0=\left \lceil \frac{d- (\delta-r)}{k+l} \right \rceil.$$  We then clearly have
$$f \leq m_0 \leq f+1.$$

{\bf Claim 1:} $\gamma_d \geq f$.

{\em Proof of Claim 1:}  Since for positive integers $s$ and $t$ we have $\lceil s/t \rceil \cdot t  \leq s+t-1$, we have
$$\lceil (l(f-1)+1)/\delta \rceil \cdot \delta \leq l(f-1)+\delta,$$
and
$$(k+l)f \leq d-\delta+(k+l)-1.$$
Therefore,
$$k(f-1) + \left \lceil (l(f-1)+1)/ \delta \right \rceil \cdot \delta \leq (k+l)(f-1) + \delta \leq d-1,$$ from which our claim follows by Lemma \ref{lem iff k,l free interval}.

{\bf Claim 2:} $\gamma_d \leq f+1$.

{\em Proof of Claim 2:}  We can easily see that
$$k(f+1) + \left \lceil (l(f+1)+1)/ \delta \right \rceil \cdot \delta > (k+l)(f+1) \geq d-\delta +k+l >d,$$
which implies our claim  by Lemma \ref{lem iff k,l free interval}.

{\bf Claim 3:} $\gamma_d \geq f+1$ if, and only if, $m_0 \geq f+1$.

{\em Proof of Claim 3:}  First note that, since $r$ is the remainder of $l f$ mod $\delta$, we have 
$$\left \lceil (lf+1)/ \delta \right \rceil \cdot \delta = lf+ \delta -r.$$
Therefore, $\gamma_d \geq f+1$ if, and only if, $$kf + lf+ \delta -r <d,$$ which is equivalent to $$f < \frac{d-(\delta -r)}{k+l};$$ since $f$ is an integer,  this is further equivalent to $f< m_0$, that is, to $f+1 \leq m_0$, as claimed.

Our result that $\gamma_d=m_0$ now follows, since if $f=m_0$, then $\gamma_d \geq f$ by Claim 1 and $\gamma_d \leq f$ by Claim 3, and if $f+1=m_0$, then $\gamma_d \geq f+1$ by Claim 3 and $\gamma_d \leq f+1$ by Claim 2. $\Box$

As a consequence of Theorem \ref{thm gamma eval}, we see the following:

\begin{cor}  \label{cor gamma lower}

For all positive integers $k$, $l$, and $d$ with $k>l$ we have \begin{eqnarray*} \gamma (\mathbb{Z}_d, \{k,l\}) & \geq & \left \lfloor \frac{d}{k+l} \right \rfloor. \end{eqnarray*}
\end{cor}

\section{Intervals Suffice}   \label{section gamma is enough}

In this section we prove Theorem \ref{gamma is enough}, that 
$$\max_{d |n} \left\{ \alpha (\mathbb{Z}_d, \{k,l\}) \cdot \frac{n}{d}  \right \} = \max_{d |n} \left\{ \gamma (\mathbb{Z}_d, \{k,l\}) \cdot \frac{n}{d}  \right \}.$$
We only need to establish that the left-hand side is less than or equal to the right-hand side, since, obviously, $$\alpha (\mathbb{Z}_d, \{k,l\}) \geq \gamma(\mathbb{Z}_d, \{k,l\}).$$  Our result will thus follow from the following:

\begin{thm}

For all positive integers $d$, $k$, and $l$ with $k>l$, there exists a divisor $c$ of $d$ for which
\begin{eqnarray*} \alpha (\mathbb{Z}_d, \{k,l\}) &\leq & \gamma (\mathbb{Z}_c, \{k,l\}) \cdot \frac{d}{c}. \end{eqnarray*}

\end{thm}

{\em Proof:}  Since $\alpha (\mathbb{Z}_d, \{k,l\})$ is the larger of $\beta (\mathbb{Z}_d, \{k,l\})$ or $\gamma (\mathbb{Z}_d, \{k,l\})$, we may assume that it equals $\beta (\mathbb{Z}_d, \{k,l\})$.  We let $\beta_d$ denote $\beta (\mathbb{Z}_d, \{k,l\})$.  

Let $A$ be a $(k,l)$-sum-free arithmetic progression in $\mathbb{Z}_d$ of size $\beta_d$, and suppose that $$A=\{a+ i \cdot b \mid i=0,1,\dots,\beta_d-1\}$$
for some elements $a$ and $b$ of $\mathbb{Z}_d$; we may assume that $\beta_d \geq 2$ (a one-element subset would be an interval) and that $g = \gcd (b,d) \geq 2$.  (We interchangeably consider $0,1,\dots,d-1$ as integers and as elements of $\mathbb{Z}_d$.)  

Let $H$ denote the subgroup of index $g$ in $\mathbb{Z}_d$.  We then have a unique element $e \in \{0,1,\dots, g -1\}$ for which $A$ is a subset of the coset $e+H$ of $H$.  We consider two cases.

When $k \not \equiv l$ mod $g$, then $\gamma_g=\gamma(\mathbb{Z}_g, \{k,l\}) \geq 1$, since (for example) $\{1\}$ is a $(k,l)$-sum-free set in $\mathbb{Z}_g$.  Therefore,         
$$\beta_d = |A| \leq |H| = d/g \leq \gamma_g \cdot d/g.$$  We thus see that $c=g$ satisfies our claim.

Assume now that $k \equiv l$ mod $g$.  In this case $ke+H=le+H$, and thus $kA$ and $lA$ are both subsets of the same coset of $H$.  Since the sets are nonempty and disjoint, we must have $|kA| < |H|$, $|lA| < |H|$, and
$$|kA|+|lA| \leq |H|.$$ 
Now $$kA=\{ka + i \cdot b \mid i=0,1,\dots, k \cdot \beta_d-k\},$$ so 
$$|kA| = \min \{|H|, \; k \cdot \beta_d-k+1\}=k \cdot \beta_d-k+1,$$ and similarly 
$$|lA| = l \cdot \beta_d-l+1.$$
Therefore, 
$$(k \cdot \beta_d-k+1)+(l \cdot \beta_d-l+1) \leq |H|=d/g,$$ from which
$$\beta_d \leq \left \lfloor \frac{d/g-2}{k+l} \right \rfloor+1.$$
Note that $\beta_d \geq 2$ implies that $$d/g-2 \geq k+l;$$ since $g \geq 2$, this then further implies that 
$$d-d/g -2 \geq k+l.$$ Therefore,
$$\beta_d \leq \left \lfloor \frac{d/g-2}{k+l} \right \rfloor +1 \leq \left \lfloor \frac{d-4}{k+l} \right \rfloor \leq \left \lfloor \frac{d}{k+l} \right \rfloor.$$  By Corollary \ref{cor gamma lower}, we thus have $\beta_d \leq \gamma_d$, which proves our claim.  $\Box$

\section{Further questions} \label{Further questions}

Having found the maximum size of $(k,l)$-sum-free sets in cyclic groups, we may turn to some other related questions.  Here we only discuss three of them; other intriguing problems, including 
\begin{itemize}
  \item the number of $(k,l)$-sum-free sets, 
  \item maximal $(k,l)$-sum-free sets (with respect to inclusion),
  \item complete $(k,l)$-sum-free sets (that is, those where $kA \cup lA =G$),
  \item maximum-size $(k,l)$-sum-free sets in subsets,
\end{itemize}
are discussed in detail in Chapter G.1.1 of the first author's book \cite{Baj:2018a}. 

\subsection{Noncyclic groups}

Clearly, if $A$ is a $(k,l)$-sum-free set in $G_1$, then $A \times G_2$ is $(k,l)$-sum-free in $G_1 \times G_2$, and thus for any abelian group of order $n$ and exponent $e(G)$ we have
\begin{eqnarray*} 
\mu (G, \{k,l\}) &\geq & \mu (\mathbb{Z}_{e(G)}, \{k,l\}) \cdot \frac{n}{e(G)}.
\end{eqnarray*}
Therefore, by Theorem \ref{mu all n k l}, 
\begin{eqnarray*} 
\mu (G, \{k,l\}) &\geq & \max_{d |e(G)} \left \{ \left \lceil \frac{d-(\delta-r)}{k+l} \right \rceil \cdot \frac{n}{d}  \right \},
\end{eqnarray*}
where $\delta=\gcd(d,k-l)$ and $r$ is the remainder of $l  \lceil (d-\delta)/(k+l) \rceil$ mod $\delta$.
We believe that equality holds.  As we mentioned in the Introduction, Green and Ruzsa proved this conjecture for the case $(k,l)=(2,1)$; see Theorem \ref{Green and Ruzsa thm} above.  As their methods were complicated and relied, in part, on a computer, we expect the general case to be challenging.  

We have the following partial result:

\begin{thm} [Bajnok, 2009; cf.~\cite{Baj:2009a}]
We have \begin{eqnarray*} 
\mu (G, \{k,l\}) &= & \mu (\mathbb{Z}_{e(G)}, \{k,l\}) \cdot \frac{n}{e(G)}
\end{eqnarray*} whenever $e(G)$ has at least one divisor $d$ that is not congruent to any integer between $1$ and $\gcd(d,k-l)$ (inclusive) mod $k+l$.
\end{thm}

In particular, for elementary abelian $p$-groups, we have:

\begin{thm}
Let $p$ be a positive prime and $r \in \mathbb{N}$.  
 If $k-l$ is divisible by $p$, then $\mu (\mathbb{Z}_p^r, \{k,l\})=0.$
 If $k-l$ is not divisible by $p$, and $p-1$ is not divisible by $k+l$, then \begin{eqnarray*} \mu (\mathbb{Z}_p^r, \{k,l\})& =& \left \lceil \frac{p-1}{k+l} \right \rceil \cdot p^{r-1}. \end{eqnarray*}
 
\end{thm}

Other cases remain open.

\subsection{Classification of maximum-size $(k,l)$-sum-free sets} 

The question that we have here is: What can one say about a $(k,l)$-sum-free subset $A$ of $G$ of maximum size $|A|=\mu (G, \{k,l\})$?

The sum-free case -- that is, when $(k,l)=(2,1)$ -- has been investigated thoroughly and is now known.  It turns out that, when the order $n$ of the group has at least one divisor that is not congruent to $1$ mod $3$, then sum-free sets of maximum size are unions of cosets that form arithmetic progressions; see the works of Diamanda and Yap in \cite{DiaYap:1969a} and Street in \cite{Str:1972a} and \cite{Str:1972b} (cf.~also Theorems 7.8 and 7.9 in \cite{WalStrWal:1972a}).  The situation is considerably less apparent, however, when all divisors of $n$ are congruent to $1$ mod $3$.  The classification was completed by Balasubramanian, Prakash, D. S. Ramana in 2016; cf.~\cite{BalPraRam:2016a}.  The general result is too complicated to present here; we just mention the example that the set 
$$\{ (n-1)/3\} \cup [(n+5)/3, (2n-5)/3] \cup \{(2n+1)/3\},$$ which is two elements short of an arithmetic progression, is sum-free in $\mathbb{Z}_n$ and has maximum size $\mu (\mathbb{Z}_n, \{2,1\})=(n-1)/3$.  (The classification of this case for cyclic groups was completed by Yap; cf.~\cite{Yap:1971a}.)

The case when $k>2$ is not known in general, but we have the following result of Plagne:

\begin{thm} [Plagne, 2002; cf.~\cite{Pla:2002a}]
Let $p$ be a positive prime, and let $k$ and $l$ be positive integers with $k>l$ and $k \geq 3$.  Suppose also that $k-l$ is not divisible by $p$.  If $A$ is a $(k,l)$-sum-free set in $\mathbb{Z}_p$ of maximum size $\left \lceil (p-1)/(k+l) \right \rceil$, then $A$ is an arithmetic progression.
\end{thm}

We are not aware of further results on the classification of $(k,l)$-sum-free sets of maximum size.

\subsection{Additive $k$-tuples}

Given a subset $A$ of $G$ and a positive integer $k$, we may ask for the cardinality $P(G,k,A)$ of the set
$$\{(a_1,\dots,a_k) \in A^k \mid a_1+ \cdots + a_k \in A\},$$ with which we can then set $P(G,k,m)$ as the minimum value of $P(G,k,A)$ among all $m$-subsets $A$ of $G$ (with $m \in \mathbb{N}$).  By definition, we have $P(G,k,m)=0$ whenever $m \leq \mu(G,\{k,1\})$, but $P(G,k,m) \geq 1$ for $\mu(G,\{k,1\}) + 1 \leq m \leq n$.

Let us consider the case of $k=2$ and the cyclic group $\mathbb{Z}_p$ of prime order $p$.  As we observed, the ``middle-third'' of the elements forms a sum-free set in $\mathbb{Z}_p$ of maximum size $\mu(\mathbb{Z}_p,\{2,1\})=\lceil (p-1)/3 \rceil$.  For $\lceil (p-1)/3 \rceil + 1 \leq m \leq p$, we may enlarge the set to 
$$A(p,m)=\{\lceil (p-m)/2 \rceil +i \mid i=0,1,\dots, m-1\}.$$ Then $A(p,m)$ is the ``middle'' $m$ elements of $\mathbb{Z}_p$, and a short calculation yields that 
$$P(\mathbb{Z}_p, 2, A(p,m))= \left \lfloor (3m-p)^2 /4 \right \rfloor.$$
Recently, Samotij and Sudakov proved that we cannot do better and that, in fact, $A(p,m)$ is essentially the only set achieving the minimum value:

\begin{thm} [Samotij and Sudakov, 2016; \cite{SamSud:2016a}, \cite{SamSud:2016b}]  \label{SamSud thm}
For every positive prime $p$ and integer $m$ with $\lceil (p-1)/3 \rceil + 1 \leq m \leq p$ we have 
\begin{eqnarray*} P(\mathbb{Z}_p, 2, m) & = &  \left \lfloor (3m-p)^2 /4 \right \rfloor.\end{eqnarray*}  Furthermore, if for some $A \subseteq G$ we have $P(\mathbb{Z}_p, 2, m)=P(\mathbb{Z}_p, 2, A)$, then there is an element $b$ of $\mathbb{Z}_p$ for which $A= b \cdot A(p,m)$.

\end{thm}

Soon after, Chervak,  Pikhurko, and  Staden generalized Theorem \ref{SamSud thm} for other values of $k$, while still remaining in cyclic groups of prime order $p$.  As they showed in \cite{ChePikSta:2017a}, the answer turns out to be more complicated, but at least in the case when $k-1$ is not divisible by $p$, the value $P(\mathbb{Z}_p, k, m)$ is still given by intervals (though there are other sets $A$ that yield the same value).  The general problem of finding $P(G,k,m)$ is largely unsolved.

\end{document}